\documentclass[12pt]{article}
 \usepackage{amsmath,amsfonts,theorem}
\input epsf
 \newcommand{\ntag}[1]{} 
\numberwithin{equation}{section}
 
 \newtheorem{prop}{Proposition}

 \newtheorem{lem}{Lemma}

 \theorembodyfont{\rmfamily}
 \newtheorem{dfn}{Definition}

 \newcommand{\qed}{\ifhmode\unskip\nobreak\fi\quad\ensuremath\square}

 \newcommand{\Oh}{\mathcal O}

 \newcommand{\al}{\alpha}
 \newcommand{\be}{\beta}
 \newcommand{\de}{\delta}

 \newcommand{\ga}{\gamma}
 
 \newcommand{\om}{\omega}

 \newcommand{\PP}{\mathbb P}
 \newcommand{\C}{\mathbb C}

 \newcommand{\R}{\mathbb R}
 \newcommand{\Z}{\mathbb Z}

\begin{document}

 \title{On families of  $BS_{can}$ - lagrangian tori in projective spaces}

  \author{N.A. Tyurin}

\date{}

\maketitle

\begin{center}
 {\em BLThP JINR (Dubna) and MSU of RT}
 \end{center}
 \bigskip

Here we discuss the following situation. Consider $\C \PP^n$ with
the standard symplectic form $\om$ such that its cohomology class
$[\om] \in H^2(\C \PP^n)$ is a generator of $H^2(\C \PP^n)$. It's
very well known, see f.e. \cite{GrHar}, that the canonical class
$K$ of the projective space equals to $-(n+1) [\om]$. Consider the
anticanonical line bundle $K^* \to \C \PP^n$ and fix a connection
$a$ on it such that the curvature form $F_a$ is proportional to
the symplectic form $\om$. Let $S \subset \C \PP^n$ is a
lagrangian torus which satisfies the following property:

\begin{dfn} A lagrangian torus is called Bohr - Sommerfeld with
respect to the anticanonical class (or $BS_{can}$ for short) if
the restriction $(K^*, a)|_S$ admits covariantly constant
sections.
\end{dfn}

It's not hard to see that the definition is correct so it doesn't
depend on the choice of $a$, see, f.e. \cite{Tyu}.

Let $S_t, t \in [0, 1]$ is a family of $BS_{can}$ - lagrangian
tori in $\C \PP^n$, the result of smooth lagrangian deformation of
a given $BS_{can}$ lagrangian torus $S_0$. Then we have the
following

\begin{lem} There exists a symplectomorphism $\phi$ of $(\C
\PP^n, \om)$ such that $\phi(S_0) = S_1$.
\end{lem}

Indeed, for every $S_t \subset \C \PP^n$ there exists the
corresponding Darboux - Weinstein neighborhood $\Oh_{DW}(S_t)
\subset \C \PP^n$ which is symplectomorphic to a neighborhood of
the zero section in $T^* S_t$. Then it follows that there is some
small $\de \in \R_+$ such that $S_l, l \in (t - \de, t+\de)$ is
contained by $\Oh_{DW} (S_t)$. Since the segment $[0,1]$ is
compact we can choose a finite covering of it by segments $(t_i -
\de_i, t_i + \de_i)$. Now one needs to show that all $BS_{can}$
lagrangian tori from our family contained by the same Darboux -
Weinstein neighborhood are symplectomorphic. To do this take a
pair $S_{\al}, S_{\be} \subset \Oh_{DW} (S_{t_i})$ of $BS_{can}$
lagrangian tori. According to the Darboux - Weinstein lemma, see
\cite{Wein}, every lagrangian torus from $\Oh_{DW} (S_t)$ is
represented by a closed 1- form on $S_t$. Since all our tori are
from the same $BS_{can}$ family it follows that both $S_{\al},
S_{\be}$ are represented by exact forms $d f_{\al},  d f_{\be}$
where $f_{\al}, f_{\be} \in C^{\infty}(S_{t_i}, \R)$ are smooth
functions on $S_{t_i}$. Of course we have in mind the
representation of our local picture inside of $T^* S_{t_i}$.
Inside of $T^* S_{t_i}$ the fact that the graphs of $d f_{\al}, d
f_{\be}$ are symplectomorphic is obvious and then we transport
this one to the Darboux - Weinstein neighborhood $\Oh_{DW}
(S_{t_i})$.

 It remains to choose a chain $S_{\al_i}, S_{\be_i}$ such that

 1. $S_{\al_i}, S_{\be_i} \in \Oh_{DW} (S_{t_i})$;

 2. $S_{\al_1} = S_0, S_{\be_{fin}} = S_1$;

 and combine the corresponding symplectomorphisms to get the desired symplectomorphism
 which maps $S_0$ to $S_1$.

As a corollary we get the following proposition

\begin{prop} Let $S_0, S_1 \subset \C \PP^n$ are $BS_{can}$ lagrangian tori
which can be joined by a family of $BS_{can}$  lagrangian tori.
Then there exists a symplectomorphism which maps $S_0$ to $S_1$.
\end{prop}

The arguments used above can be extended to the case when two
given $BS_{can}$ lagrangian tori can be joined by a family which
can contain  non $BS_{can}$ lagrangian tori.

Consider lagrangian family $S_t, t \in [0,1]$ such that almost all
$S_t$ are $BS_{can}$ except some small subsegment $S_t, t \in (t_0
- \de, t_0 + \de)$ of tori contained by the Darboux - Winstein
neighborhood of certain $S_{t_0}$. For this situation we can
extend the arguments used above modulo certain correction which
comes with certain characteristic class discussed in \cite{Tyu}.
Recall that a $BS_{can}$  lagrangian submanifold $S$ of a monotone
simply connected symplectic manifold admits an integer class $m_S
\in H^1(S, \Z)$ which was called the universal Maslov class.

The extension of Lemma 1 looks as follows

\begin{lem} Let $S_0, S_1 \subset \C \PP^n$ are $BS_{can}$
lagrangian tori which can be joined by some lagrangian family
$S_t, t \in [0,1]$ . Let $S_{t_0} \subset \C \PP^n$ is a
lagrangian torus from the family such that both $S_0, S_1$ are
contained by its Darboux - Weinstein neighborhood $\Oh_{DB}(S)$.
Then $S_0$ and $S_1$ are symplectomorphic if and only if $m_{S_0}
= m_{S_1}$.
\end{lem}

Note, that since $m_{S_i} \in H^1(S_i, \Z)$ one could not compare
these classes in general situation. But here we have some
lagrangian family, joined $S_0$ and $S_1$, and thus one has an
identificaion of the groups.

Consider the Darboux - Weinstein neighborhood $\Oh_{DW}(S_{t_0})$
and transport $S_0$ and $S_1$ to the corresponding neighborhood of
the zero section inside $T^* S_{t_0}$. We denote the images by the
same symbols. Over $S_{t_0}$ the (images of) lagrangian
submanifolds are presented by two closed 1 -forms, $\al_0$ and
$\al_1$. $S_0$ and $S_1$ are symplectomorphic if and only if the
forms represent the same cohomology class in $H^1(S_{t_0}, \R)$.
Remark that $S_0$ and $S_1$ have the same Maslov indecis being in
the same family. It follows that every pair of loops $\ga_0,
\ga_1$, $\ga_i \subset S_i$ and $[\ga_0] = [\ga_1]$ under the
identification mentioned above, is joined by tube of zero
symplectic area if and only if $m_{S_0} = m_{S_1}$. Thus $m_{S_0}
= m_{S_1}$ if and only if for any pair of cohomologically
equivalent loops $\ga_0, \ga_1, \ga_i \subset S_i$
$$
\int_{\ga_0} \rho = \int_{\ga_1} \rho
$$
by the Stokes formula, where $\rho$ is the canonical 1 - form on
$T^* S_{t_0}$. Further, by the definition of the canonical 1 -
form on $T^* S_{t_0}$ we have that for any loop $\ga \subset
S_{t_0}$
$$
\int_{\ga} \al_0 = \int_{\ga} \al_1.
$$
It is possible if and only if $[\al_0] = [\al_1] \in H^1(S_{t_0},
\R)$. This means that
$$
\al_0 - \al_1 = d f,
$$
where $f$ is a smooth function on $S_{t_0}$ whose hamiltonian
vector field generates the desired symplectomorphism which maps
$S_0$ to $S_1$.

Now we can extend slightly the statement of Proposition 1.

Let $S_t, t \in [0,1]$ is a lagrangian family with {\it small} non
$BS_{can}$ pieces, so every non $BS_{can}$ piece is covered by the
Darboux - Weinstein neighborhood of a member. Then one has

\begin{prop} Let $S_0, S_1 \subset \C \PP^n$ be $BS_{can}$
lagrangian tori in the projective space. If there exists a
lagrangian family with small non $BS_{can}$ pieces and such that
all $BS_{can}$ members have the same universal Maslov class, then
$S_0$ is symplectomorphic to $S_1$.
\end{prop}

The simple remarks, presented above, could be exploited to study
lagrangian families which join lagrangian tori of Clifford and
Chekanov types in $\C \PP^2$ in the representation, given in
\cite{Aur}, and as a corollary one would get that a $BS_{can}$
torus of Chekanov type is not monotone.

\end{document}